\documentclass[12pt]{amsart}

\usepackage[utf8]{inputenc}

\usepackage{caption}
\usepackage{tikz}

\usepackage{subcaption}
\usepackage{setspace}
\usepackage{eqnarray,amsmath}
\usepackage{enumerate}
\usepackage[shortlabels]{enumitem}
\usepackage{amsthm, amsfonts,amssymb, stmaryrd,yfonts,pxfonts,pifont, bbm,tikz-cd}
\usepackage{caption}
\usepackage{enumitem}

\newcommand{\N}{\mathbb{N}}

\newcommand{\R}{\mathbb{R}}
\newcommand{\Z}{\mathbb{Z}}

\usepackage{graphicx}

\usetikzlibrary{trees}

 \usepackage{hyperref}

\usepackage{latexsym}
\usepackage{array}
\usepackage{xcolor}

\newtheorem{theorem}{Theorem}

\newtheorem{cor}{Corollary}

\newtheorem{example}{Example}
\newtheorem{remark}{Remark}

\theoremstyle{definition}

\begin{document}
\title{Dynamics of Enveloping Semigroup of Flows}
\author{Sushmita Yadav and Puneet Sharma}
\address{Department of Mathematics, I.I.T. Jodhpur, N.H. 65, Nagaur Road, Karwar, Jodhpur 342037, INDIA}%
\email{yadav.34@iitj.ac.in, puneet@iitj.ac.in}%

\begin{abstract}
In this article, we relate the dynamics of a flow $(X, T)$ with the dynamics of the induced flow $(E(X), T)$ where $E(X)$ is the enveloping semigroup of flow $(X, T)$. We establish that a flow $(X, T)$ is distal if and only if the induced flow is also distal. We prove that while the induced flow on the enveloping semigroup of an equicontinuous flow is equicontinuous, the converse holds when $(X, T)$ is point transitive. We prove that for any flow $(X,T)$, $(E(X),T)$ is isomorphic to $(E(E(X)),T)$. We show that if a distal flow contains a minimal sensitive subsystem, the induced flow on the enveloping semigroup is also sensitive. We relate various forms of rigidity for a flow with rigidity for the induced flow. We also establish that for any proximal flow $(X, T)$, if T is abelian then the induced flow $(E(X),T)$ is also proximal.
 \end{abstract}
\maketitle
\section*{Introduction}
By a \emph{topological flow} we mean $(X, T,\pi)$, where $X$ is a compact Hausdorff topological space without any isolated points, $T$ is a topological group, and $\pi: T \times X \rightarrow X$  is a continuous function satisfying $\pi(st,x)=\pi(s,\pi(t,x))$. For the sake of notation, $\pi(t,x)$ is denoted as $t.x$ or $tx$ and describes the action of the group $T$ on $X$. For each $t\in T$, $\pi^t:X\rightarrow X$ defined as $\pi^t(x)=tx$ is a homeomorphism on $X$ and the set  $\{\pi^{t}\ |\ t\in T\}$ forms a group of homeomorphisms of $X$ under the composition. 
% In particular, if $T = \R$ (or $\Z$), we call the flow a \emph{continuous (or discrete)}  dynamical system. In the case of a discrete dynamical system, we will denote the flow $(X, T,\pi)$ as $(X, f)$, where $\pi(n,x)=f^{n}(x)$.
If the map $\pi$ is understood, we abbreviate the flow as $(X, T)$. In case the group $T$ is abelian, we say the flow $(X, T)$ is abelian.
For any flow $(X,T)$, the \emph{Enveloping semigroup} of $(X, T)$ is defined as $E(X,T)= \overline{\{\pi^t:X \rightarrow X | t\in T\}}$ (where $X^X$ is equipped with product topology). In case the acting group $T$ is understood, we denote the enveloping semigroup  by $E(X)$. It may be noted that for any $t\in T,~ p\in E(X)$, $\hat{\pi}(t,p)=\pi^t\circ p~ (=t.p)$ defines an action of the group $T$ on $E(X)$. Consequently $(E(X),T)$ forms a \emph{topological flow}.\\

For any topological flow $(X,T)$, the set $Tx=\{t.x: t\in T\}$ is called \emph{orbit} of $x\in X$ and denoted by $\mathcal{O}(x)$. A subset $A\subseteq X$ is said to be \emph{invariant} if $TA=\{ta:t\in T,a\in A\}\subseteq A$. If $A$ is invariant, then $(A, T)$ forms a flow and is said to be a \emph{subflow} of $(X, T)$. Any flow $(X,T)$ is point transitive if there exists $x\in X$ such that $\mathcal{O}(x)$ forms a dense subset of $X$. Any such $x$ is called a \emph{transitive point} for $(X,T)$. By $Trans(X)$, we denote the set of transitive points for $(X,T)$. A flow $(X, T)$ is called \emph{topologically transitive} if for every pair of non-empty open sets $U, V$ in $X$, there exists $t\in T$ such that $tU\cap V\neq \emptyset$. A subset $M\subseteq X$ is said to be \emph{minimal} if it is a non-empty, closed-invariant subset of $X$ with no proper subset having these properties. A subset $S \subseteq T$ is said to be \emph{syndetic} if there exists a compact subset $K\subset T$ such that $KS=T$. A subset $S$ of $T$ is called \emph{thickly syndetic} if for every compact subset $K \subset T$, there is a syndetic set $S_K \subset T$ such that $KS_K\subset S$.
\\

\noindent Let $X$ be a set and $\mathcal{U}_X$ be the collection of subsets of $X\times X$ satisfying the following:
\begin{itemize}
\item $\emptyset\notin \mathcal{U}_X$ and $X\times X\in \mathcal{U}_X$.
\item For any $U\in \mathcal{U}_X$, $\Delta\subset U$ (where $\Delta$ is the diagonal subset of $X\times X$).
\item If $U_1,U_2\in \mathcal{U}_X$ then $U_1\cap U_2\in \mathcal{U}_X$.
\item If $U\in \mathcal{U}_X$ and $V\subseteq X\times X$ such that $U\subseteq V$ then $V\in \mathcal{U}_X$.
\item For any $U\in \mathcal{U}_X$ there exists $V\in \mathcal{U}_X$ such that $U=V\circ V$ (where $V_1\circ V_2=\{(x,y)\in X\times X : \exists z\in X ~such~ that~(x,z)\in V_1~and~ (z,y)\in V_2\}$).
\end{itemize}
Then $(X,\mathcal{U}_X)$ is said to be \emph{uniform space} and $\mathcal{U}_X$ is called \emph{uniformity} on $X$. Any element $U\in \mathcal{U}_X$ is referred as an entourage. It is known that every compact Hausdorff topological space $(X,\tau)$ is equipped with unique uniformity (compatible with $\tau$). For any $x\in X,~U\in\mathcal{U}_X$, let $U[x]=\{y\in X:(x,y)\in U\}$. Then, the collection $\{U[x]:U\in \mathcal{U}_X\}$ forms a neighbourhood system for $x$ (in $X$). \\

A pair $(x,y)\in X\times X$ is said to be \emph{proximal} in $(X,T)$ if for every $U\in \mathcal{U}_X$, there exists $t\in T$ such that $(tx,ty)\in U$. The pair $(x,y)$ is called \emph{syndetically proximal} if for every $U\in \mathcal{U}_X$, the set $\{t\in T: (tx,ty)\in U\}$ is a syndetic subset of $T$. A flow is said to be proximal (syndetically proximal) if every pair in $X\times X$ is proximal (syndetically proximal). A pair $(x,y)\in X\times X$ is said to be \emph{distal} in $(X,T)$ if it is not a proximal pair. Any flow $(X, T)$ is called distal if every non-trivial pair in $X\times X$ is distal. A flow $(X,T)$ is \emph{syndetically distal} if no non-trivial pair in $X\times X$ is syndetically proximal. A flow $(X,T)$ is said to be \emph{equicontinuous} at $x\in X$ if for each $V\in \mathcal{U}_X$ there exists $U\in \mathcal{U}_X$ such that $(x,y)\in U$ implies $(tx, ty)\in V$ for all $t\in T$. A flow $(X, T)$ is said to be equicontinuous if it is equicontinuous at each $x\in X$. Let $Eq(X)$ denote the set of points of equicontinuity for $(X,T)$. If the flow $(X, T)$ has a dense set of equicontinuity points, then $(X, T)$ is referred as an \emph{almost equicontinuous} flow. Any flow $(X, T)$ is called \emph{hereditarily almost equicontinuous} if every closed sub-flow of $(X, T)$ is almost equicontinuous. \\

A flow $(X,T)$ is said to be \emph{sensitive} at $x\in X$ if there exists $V\in\mathcal{U}_X$ such that for every $U\in \mathcal{U}_X$, there exists $y\in U[x]$ and $t\in T$ such that $(tx,ty)\notin V$. A flow $(X,T)$ is said to be sensitive if there exists  $V\in\mathcal{U}_X$ such that for every $x\in X$ and $U\in \mathcal{U}_X$ there exists $y\in U[x]$ and $t\in T$ such that $(tx,ty)\notin V$. A flow $(X,T)$ is said to be \emph{weakly rigid} if for every $n\in\N$, if $x_1,x_2,\ldots,x_n\in X$ and $U\in \mathcal{U}_X$ then there exists $t\in T$ such that $(tx_i,x_i)\in U$ for all $i=1,2,\ldots,n$. The flow $(X, T)$ is called \emph{uniformly rigid} if there exists a net $(t_{\alpha})$ such that $(t_{\alpha})$ converges uniformly to $e$. Any continuous surjective map $\phi:(X,T)\rightarrow (Y,T)$ is said to be a 
homomorphism or factor map if $\phi(tx)=t\phi(x)$ for all $t\in T$. In this case, the flow $(Y,T)$ is called a factor of $(X,T)$. Further, if $\phi$ is bijective, the flows $(X,T)$ and $(Y,T)$ are said to be isomorphic to each other (denoted as $(X,T)\cong (Y,T)$). Any map $f: X\rightarrow Y$ is said to be feeble open if for any non-empty open set $U\subseteq X$, $int(f(U))\neq \emptyset$.\\

The notion of enveloping semigroup was introduced by R. Ellis in $1960$ (\cite{ellis}). He established that for any flow $(X,T)$, proximality is an equivalence relation if and only if $E(X)$ contains a unique minimal right ideal.  In \cite{en}, the authors established that $(X, T)$ is weakly almost periodic if and only if every element in $E(X)$ is continuous. They also proved that if the topological flow $(X, T)$ is minimal and every element in $E(X)$ is continuous then the flow $(X, T)$ is equicontinuous. In \cite{el}, the author investigated the enveloping semigroup of flow generated by projective space of a finite-dimensional vector space, when the acting group is a subgroup of the general linear group $GL(V)$. In \cite{ko}, the author fully characterized the enveloping semigroup of strongly non-chaotic systems which contains periodic points of period $2^k$ (for every $k\in\N$). Although computation of the enveloping semigroup is often cumbersome, the enveloping semigroup for some of the cases has been computed explicitly  \cite{gla,fur,el}. In \cite{gmu}, authors gave necessary and sufficient criteria for the metrizability of the enveloping semigroup of a flow. In \cite{glas}, Glasner studied tame dynamical systems and established that for minimal distal but not equicontinuous systems, the canonical map from the enveloping operator semigroup onto the Ellis semi-group is never an isomorphism. He also gave a complete characterization of minimal systems whose enveloping semigroup is metrizable. In \cite{qu}, the authors investigated the topological flow generated by discrete group action on a countable compact metric space and established that the equicontinuity and distality of such flows are equivalent.\\

 % It is known that the enveloping semigroup of a flow also constitutes a flow, a natural question arises regarding the dynamic relationship between these two structures.

In this paper, we relate the dynamics of a flow $(X,T)$ with the dynamics of the induced flow $(E(X),T)$. We establish that a flow $(X,T)$ is distal if and only if the flow $(E(X),T)$ is distal. We prove that if a flow is equicontinuous, the induced flow $(E(X),T)$ is also equicontinuous. Further, if $(X,T)$ is point transitive then equicontinuity of $(E(X),T)$ guarantees equicontinuity of the underlying system. We prove that for any flow $(X, T)$, the flow $(E(X),T)$ is isomorphic to $(E(E(X)),T)$. We also relate sensitivity of a flow with sensitivity of its enveloping semigroup. We prove that a distal flow containing a minimal sensitive subsystem has a sensitive enveloping semigroup. Finally, we relate various forms of rigidity for a flow with similar notions for the induced flow.

\section*{Main Results}

\begin{theorem}\label{dis}
$(X, T)$ is distal if and only if $(E(X),T)$ is distal.
\end{theorem}

\begin{proof} 
Let $(X,T)$ be distal and $(p,q)$ be a proximal pair for  $(E(X), T)$. As $(p,q)$ is proximal, there exists a net $(t_{\alpha})$ in $T$ such that $\lim ~t_{\alpha}p = \lim ~t_{\alpha}q$ and hence for any $x\in X$, $\lim ~ t_{\alpha}px = \lim ~ t_{\alpha}qx$. Consequently, the pair $(px,qx)$ is proximal for $(X,T)$ and thus $px=qx$. As the proof holds for any $x\in X$, $p=q$ and the proof of the forward part is complete.  

Conversely, if $(x,y)$ is proximal for $(X,T)$ then there exists a minimal ideal $K \subseteq E(X)$ such that $px=py$ for any $p\in K$. As $E(X)$ is distal, $e\in K$ (as $e$ is the only idempotent). Thus, we have $x=y$ and the proof is complete. 
\end{proof}

\begin{theorem}\label{group}
Given a flow $(X, T)$, $E(X)$ forms a group of homeomorphism $\implies$ $E(E(X))$ forms a group of homeomorphism.
\end{theorem}
\begin{proof}
As $(X, T)$ is distal if and only if $(E(X), T)$ is distal (Theorem~\ref{dis}), we have $E(X)$ is a group if and only if $E(E(X))$ forms a group. Let $\hat{g}\in E(E(X))$, $\xi\in E(X)$ and $(p_\alpha)$ be a net such that $p_\alpha\rightarrow \xi$ in $E(X)$. As $\hat{g}\in E(E(X))$, there exists $g\in E(X)$ such that $\hat{g}(p)=g\circ p$ for all $p\in E(X)$. As $E(X)$ is a group of homeomorphisms, $\hat{g}(p_\alpha)=g\circ p_\alpha\rightarrow g\circ \xi=\hat{g}(\xi)$. Hence, $\hat{g}$ is continuous (and hence a homeomorphism), and the proof is complete.
\end{proof}

\begin{theorem}
 $(X,T)$ is equicontinuous then $(E(X),T)$ is equicontinuous. If the flow $(X, T)$ is point-transitive, then the converse also holds.
\end{theorem}

\begin{proof}
It may be noted that if $(X, T)$ is equicontinuous, then $E(X)$ forms a group of homeomorphisms and hence $E(E(X))$ forms a group of homeomorphisms (Theorem~\ref{group}). Consequently, $(E(X), T)$ forms an equicontinuous flow. Conversely, if $(X, T)$ is point-transitive then for any transitive point $x_0\in X$, the map $\phi(p)=px_0$ defines a factor map from $(E(X),T)$ onto $(X,T)$. As a factor of equicontinuous flow is equicontinuous~\cite{a}, $(X, T)$ is equicontinuous.
\end{proof}

\begin{theorem}\label{t3}
If $(X, T)$ is metrizable, distal, and hereditarily almost equicontinuous, then $(E(X), T)$ forms an equicontinuous flow.
\end{theorem}
\begin{proof}
Since $(X, T)$ is hereditarily almost equicontinuous, $(E(X), T)$ is metrizable hereditarily almost equicontinuous \cite[Corollary 10.2]{g}. Also, as $(X, T)$ is distal, $(E(X), T)$ is minimal. Further, as any minimal flow is either equicontinuous or sensitive, $(E(X), T)$ forms an equicontinuous flow, and the proof is complete.
\end{proof}

\begin{remark}
The above results relate distality and equicontinuity for the flows $(X, T)$ and $(E(X), T)$. The results establish that while distality is equivalent for the two systems, equicontinuity of $(X, T)$ ensures equicontinuity for $(E(X), T)$. Also, Theorem \ref{t3} establishes that if a metric flow $(X, T)$ is distal and hereditarily almost equicontinuous, then $(E(X), T)$ forms an equicontinuous flow. It may be noted that equivalence of equicontinuity for the two systems $(X, T)$ and $(E(X), T)$ holds strictly under point transitivity of $(X, T)$ and the result need not hold if the assumption is dropped (Example \ref{circ}). It may be noted that as the natural correspondence $\pi^t\Leftrightarrow \hat{\pi}^t$ provides a conjugacy between $(E(X),T)$ and $(E(E(X)),T)$, the systems $(E(X),T)$ and $(E(E(X)),T)$ are conjugate (to each other). We now establish our claims below. %In \cite{NM}, the authors computed $E(X)$ and $E(E(X))$ of the discrete flow $(X,h)$ where $X=\{(r,\theta): 1\leq r \leq 2, \theta \in [0,1)\}$ and $h: X \rightarrow X$ be defined as $h(r,\theta)=(1+(r-1)^2,\theta+\alpha)$. Further, they established that $(E(X),T) \ncong (E(E(X)),T)$. 
\end{remark}

\begin{example}\label{circ}
Let $X =\{(r,\theta) : r \in \{2-\frac{1}{2^n} : n\in N\}\cup \{1,2\}, 0\leq \theta \leq 1\}$ and $f:X\rightarrow X$ is defined as $f(r,\theta)=(r,\theta + r)$ (represented in Figure ~\ref{fig:1}).
 \begin{figure}[hbt!]
    \centering
    \includegraphics[width=0.7\linewidth]{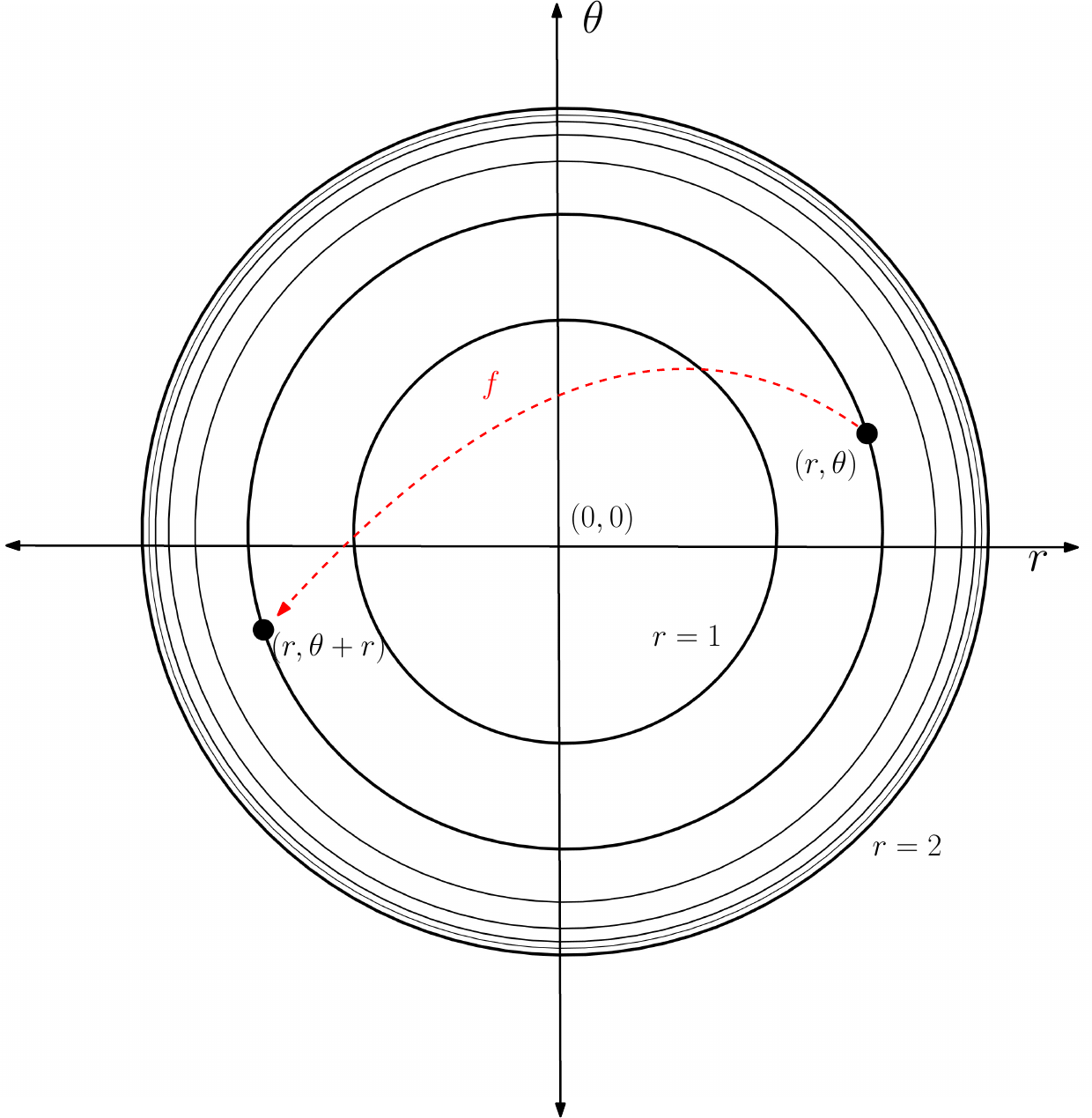}
    \caption{Dynamical Plot for Example $1$}
    \label{fig:1}
\end{figure} 
Firstly, note that the map $f$ acts as an identity on the innermost circle ($r=1$). For $r=\frac{3}{2}$, while the even iterates of $f$ act as identity, the odd iterates rotate the points by $\theta=\frac{1}{2}$ (with points on the innermost circle unchanged under $f$ in either case). For the $k$-th stage, $\mathbb{N}$ can be partitioned into $2^k$ sets $S_{k0},S_{k1},\ldots,S_{k2^k-1}$ where $S_{ki}=\{n\in\mathbb{N}: n\cong i~~ mod (2^k)\}$ and the membership of any $r$ (in one of $S_{ki}$'s) completely determines the action of $f^r$ (and hence it's iterates) on the first $k$ circles. 
\begin{figure}[hbt!]
    \centering
    \includegraphics[width=0.9\linewidth]{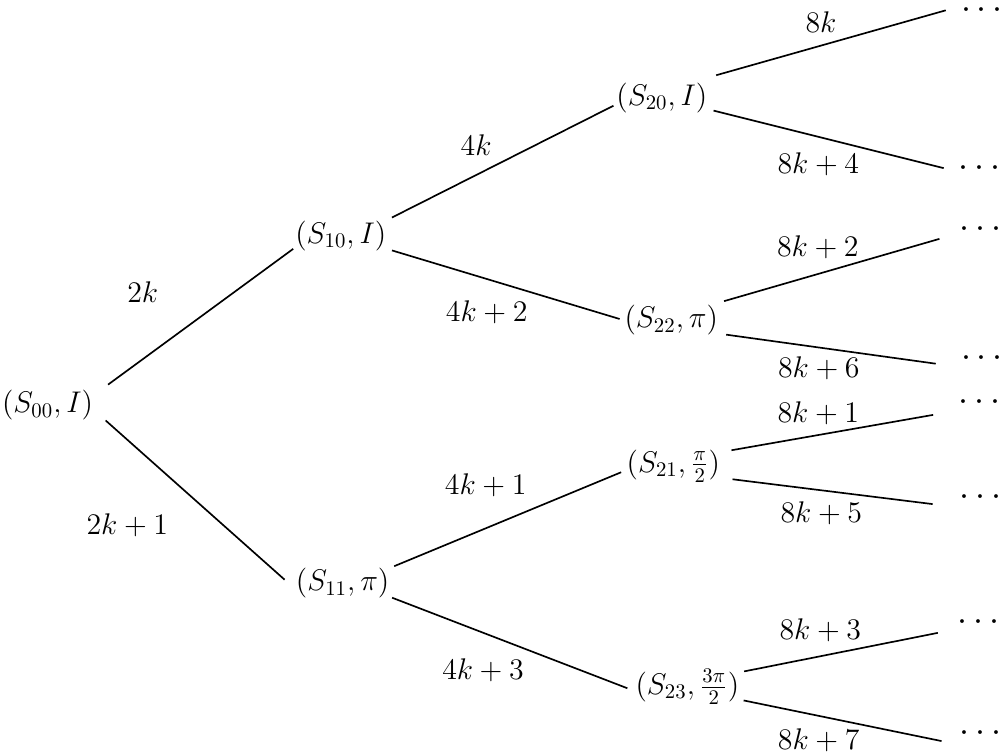}
    \caption{Computation of $E(X,f)$ for Example \ref{circ}}
    \label{fig:2}
\end{figure}
It can be seen that elements of $E(X)$ are nothing but paths along the Figure~\ref{fig:2} (where the action on the $k$-th circle is nothing but the $k$-th entry along the path). It is clear that as $(X,f)$ is metrizable and hereditarily almost equicontinuous, $E(X)$ is metrizable and hereditarily almost equicontinuous. One can define the metric $d'$ on $E(X)$ as $$d'(g_1,g_2)=sup_{\theta\in[0,2\pi)}\lvert g_1(2,\theta)-g_2(2,\theta)\rvert+\sum_{n=1}^{\infty}\frac{1}{2^n}sup_{\theta\in[0,2\pi)]} \lvert g_1(r_n,\theta)-g_2(r_n,\theta)\rvert,$$ where $r_n=2-\frac{1}{2^n}.$ One can verify that metric $d'$ generates the product topology on $E(X)$. As $(E(X), T)$ is minimal and hereditarily almost equicontinuous, $(E(X), T)$ is equicontinuous. Thus, the example provides an instance when $(E(X), T)$ is equicontinuous but $(X, T)$ is not equicontinuous. 
\end{example}

\begin{theorem}\label{iso}
For any flow $(X,T)$, $(E(X),T)\cong (E(E(X)),T)$.
\end{theorem}
\begin{proof}  Let $E(X)=\{\overline{\pi^{t} : t\in T}\}$ and $E(E(X))=\{\overline{\hat{\pi}^{t} : t\in T}\}$, where $\pi^{t}(x)=tx$ and  $\hat{\pi}^{t}(g)=tg=\pi^{t}\circ g$. 
First, we first prove that $\pi^{t_\alpha}\rightarrow g$ in E(X) if and only if $\hat{\pi}^{t_{\alpha}}\rightarrow \hat{g}$ in E(E(X)) where $\hat{g}(p)=g\circ p$ for all $p\in E(X)$. Firstly, assume that $\pi^{t_\alpha}\rightarrow g$ and $p\in E(X)$. For any $p$ in $E(X)$, recall that $\hat{\pi}^{t_{\alpha}}(p)=\pi^{t_{\alpha}}\circ p$. As $\pi^{t_{\alpha}}\circ p(x)\rightarrow g(p(x))$ for all $x\in X$, we have $\pi^{t_{\alpha}}\circ p\rightarrow g\circ p$ for all $p\in E(X)$ and hence $\hat{\pi}^{t_{\alpha}}(p)\rightarrow \hat{g}(p)$. Conversely, if $\hat{\pi}^{t_{\alpha}}\rightarrow \hat{g}$ then $\hat{\pi}^{t_{\alpha}}(p)\rightarrow \hat{g}(p)$ for all $p\in E(X)$ which implies that  $\pi^{t_{\alpha}}\circ p\rightarrow g\circ p$ for all $p\in E(X)$. In particular, take $p=e$, then we have $\pi^{t_{\alpha}}\rightarrow g$.\\ 
Recall that T acts $E(E(X))$, through the map $\hat{\hat{\pi}}^t:E(E(X))\rightarrow E(E(X))$ defined by $\hat{\hat{\pi}}^t(\hat g)=\hat{\pi}^t\circ \hat{g}$. Then, the map $G: (E(X),T) \rightarrow (E(E(X)),T)$ defined as $G(\lim ~ \pi^{t_{\alpha}})=\lim ~ \hat{\pi}^{t_{\alpha}}$ is a well defined bijective continuous map such that $\hat{\hat{\pi}}^t\circ G=G \circ \hat{\pi}^t$ $\forall~~t\in T$. Futher, for any $g=\lim \pi^{t_\alpha}$, $G(\hat{\pi}^t(\lim \pi^{t_\alpha}))=G(\pi^t\circ \lim\pi^{t_\alpha})=G(\lim(\pi^{t\ast t_\alpha}))= \lim \hat{\pi}^{t\ast t_\alpha}$ and $\hat{\hat{\pi}}^t(G(\lim \pi^{t_\alpha}))=\hat{\hat{\pi}}^t(\lim \hat{\pi}^{t_\alpha})=\hat{\pi}^t\circ \lim \hat{\pi}^{t_\alpha}= \lim \hat {\pi}^t \circ \hat {\pi}^{t_\alpha}$. It is easy to see that $\hat{\pi}^{t\ast t_\alpha}(p)=\hat {\pi}^t \circ \hat {\pi}^{t_\alpha}(p)$, $\forall ~ p\in E(X)$. Hence $G$ is an isomorphism and $(E(X),T)\cong (E(E(X)),T)$.
\end{proof}

\begin{remark}
The above result establishes that for any flow $(X,T)$, $(E(X),T)$ and $(E(E(X)),T)$ are conjugate to each other. The proof establishes the natural correspondence as the conjugacy between the two systems and thus provides the change of coordinates under which one system can be visualized as the other. Consequently, the dynamical behavior of $(E(X),T)$ is replicated by its induced system, and thus, the dynamical complexity does not magnify (when compared with further induced systems). We now illustrate the statement with the help of an example.  

\end{remark}
\begin{example}\label{annulus} 
Let $X=\{(r,\theta): 1\leq r \leq 2, \theta \in [0,1)\}$ (annular region with total angle scaled to $1$), $\alpha$ be an irrational and $h: X \rightarrow X$ be defined as $h(r,\theta)=(1+(r-1)^2,\theta+\alpha)$ (illustrated in Fig.\ref{fig:3}).
\begin{figure}[hbt!]
    \centering
    \includegraphics[width=0.5\linewidth]{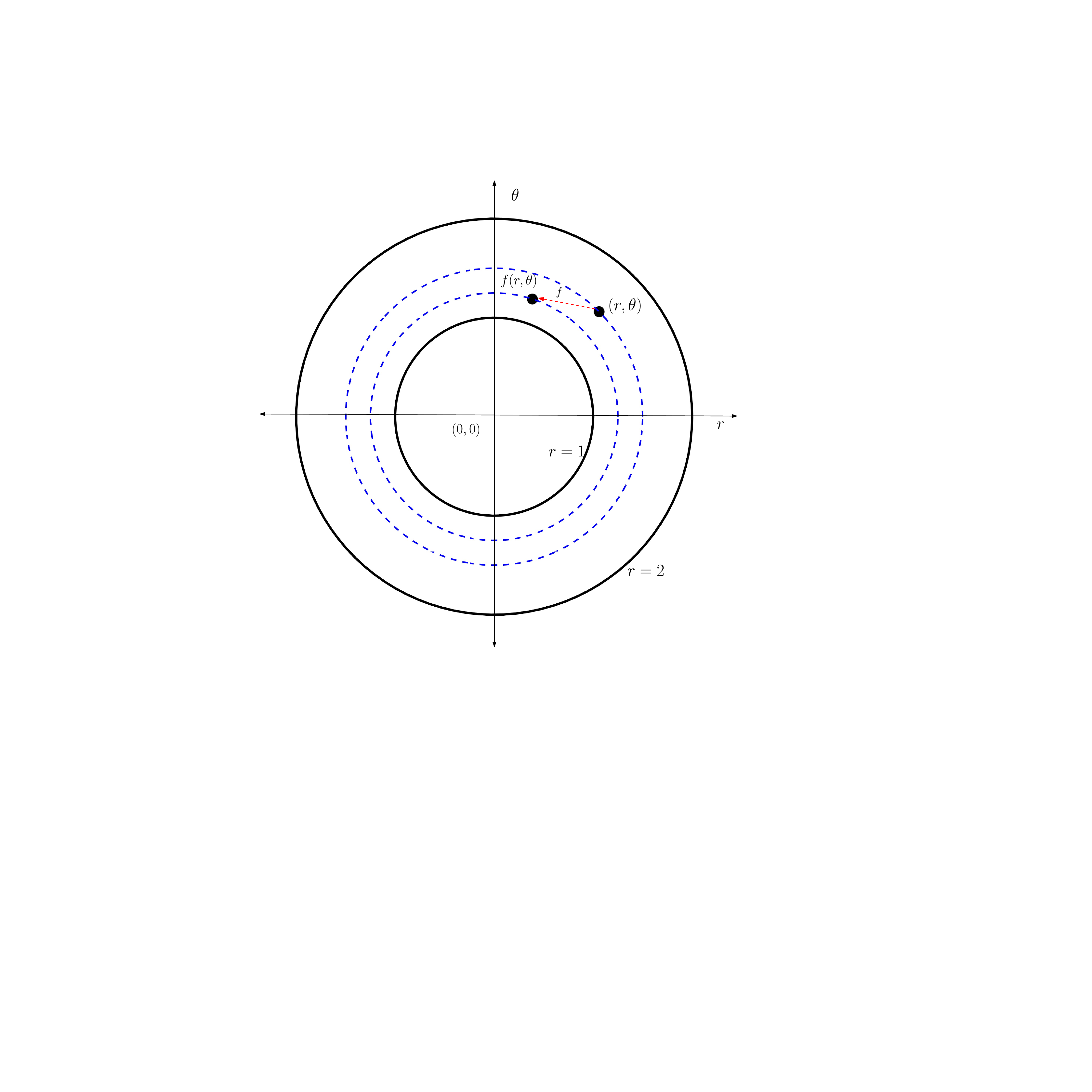}
    \caption{Dynamical Plot for Example \ref{annulus}}
    \label{fig:3}
\end{figure}
As $\alpha$ is irrational, for any $\beta\in [0,1)$, there exists a sequence $(n_k)$ of natural numbers such that $n_k\alpha~~(mod 1)$ converges to $\beta$. Consequently, $h^{n_k}\rightarrow h_{1\beta}$ and $h^{-n_k}\rightarrow h_{2\beta}$ where\\

$h_{1\beta}(r,\theta) = \left\{
     \begin{array}{lr}
       {(1,\theta+\beta)} &  :r\in [1,2) \\
       {(2,\theta+\beta)} & : r=2 \\
      
     \end{array}\right\}$ and $h_{2\beta}(r,\theta) = \left\{
     \begin{array}{lr}
       {(1,\theta-\beta)} & : r=1 \\
       {(2,\theta-\beta)} &  :r\in (1,2] \\
      
     \end{array}
   \right\}$
\vspace{0.3cm}\\

Hence, E(X)=$\{h^n:n\in \Z\}\cup \{h_{1\beta}:\beta\in[0,1)\}\cup\{h_{2\beta}:\beta\in[0,1)\}$ and the operations in $E(X)$ can be captured through the following relations: \\

\noindent
$hh_{1\phi}= h_{1(\phi+\alpha)},~~ hh_{2\phi}=h_{2(\phi-\alpha)},~~ h_{1\phi}h_{2\phi}=h_{20},~~ h_{2\phi}h_{1\phi}=h_{10},\\ h_{1\phi}h_{2\chi}=h_{2(\chi-\phi)},~~ h_{2\chi}h_{1\phi}=h_{1(\phi-\chi)},~~ h_{2\phi}h_{2\chi}=h_{2(\phi+\chi)},~~ h_{1\phi}h_{1\chi}=h_{1(\chi+\phi)}$.\\
%Note: $h_{10}$ and $h_{20}$ are only idempotents in E(X).\\

To compute $E(E(X))$, note that as $\Z$ acts on $E(X)$, each $n\in \Z$ determines map $\hat {h}^n:E(X)\rightarrow E(X)$ defined as $\hat{h}^n(g)=h^n\circ g$ and $E(E(X))=\overline{\{\hat{h}^n:n\in \Z\}}$. Once again, as for each $\eta\in [0,1)$, there exists a sequence $(n_k)$ of natural numbers such that $n_k\alpha~~(mod 1)$ converges to $\eta$, we have the following identities:\\

$\hat{h}^{n_k}(h^m)=h^{n_k}\circ h^m\rightarrow h_{1(\eta+ m\alpha)}$\\
$\hat{h}^{n_k}(h_{1\phi})=h^{n_k}\circ h_{1\phi}=h_{1(\phi+ n_k\alpha)}\rightarrow h_{1(\phi+\eta)}$\\
$\hat{h}^{n_k}(h_{2\phi})=h^{n_k}\circ h_{2\phi}=h_{2(\phi-n_k\alpha)}\rightarrow h_{2(\phi-\eta)}$\\

Thus, for each $\eta \in [0,1)$, there exists a sequence $(n_k) \in \N$ such that $\hat{h}^{n_k}\rightarrow H_{1\eta}$ and $\hat{h}^{-n_k}\rightarrow H_{2\eta}$ where $H_{1\eta}: E(X)\rightarrow E(X)$ and $H_{2\eta}:E(X)\rightarrow E(X)$ are defined as:\\

$H_{1\eta}(g) = \left\{
     \begin{array}{lr}
       {h_{1(\eta+k\alpha)}} &  : g=h^k \\
       {h_{1(\phi+\eta)}} & :g=h_{1\phi} \\
       {h_{2(\phi-\eta)}} & :g=h_{2\phi}\\
     \end{array}
   \right.$ and  $H_{2\eta}(g) = \left\{
      \begin{array}{lr}
       {h_{2(\eta-2\pi k\alpha)}} &  : g=h^k \\
       {h_{1(\phi-\eta)}} & :g=h_{1\phi} \\
       {h_{2(\phi+\eta)}} & :g=h_{2\phi}\\
     \end{array}
   \right.$
\vspace{0.3cm}\\
\noindent
Hence, $E(E(X))=\{\hat{h}^n,n\in \Z\}\cup \{H_{1\eta}:\eta\in[0,1)\} \cup \{H_{2\eta}:\eta\in[0,1)\}$. Further, as the diagram 

$$\begin{tikzcd} E(X)\arrow[r,"\hat{h}"]\arrow[d,"G"] &E(X)\arrow[d,"G"]\\ E(E(X))\arrow[r,"\hat{\hat{h}}"] &E(E(X)) \end{tikzcd}$$

commutes where $G(\lim ~ \pi^{t_{\alpha}})=\lim ~ \hat{\pi}^{t_{\alpha}}$ and $\hat{\hat{h}}^n(g)=\hat{h}^n\circ \hat{g}$, we have $(E(X),T)\cong (E(E(X)),T)$ and the proof is complete.
\end{example}

\noindent Before we move ahead, we provide a known  criteria to establish weak rigidity for general flows.

\begin{theorem}\cite{NM2}\label{niso}
For any flow $(X, T)$, the following are equivalent: 
\begin{enumerate}
\item (X, T) is weakly rigid.
\item The identity $e$ is not isolated in E(X).
\end{enumerate}
\end{theorem}

\begin{cor}\label{wr}
If $(E(X), T)$ is sensitive, then $(X, T)$ is weakly rigid.
\end{cor}

\begin{proof}
Since $(E(X), T)$ is sensitive, $E(X)$ does not have any isolated point and hence $e\in E(X)$ is not isolated. Consequently, $(X, T)$ is a weakly rigid flow (theorem~\ref{niso}).
\end{proof}

\begin{cor}
Enveloping semigroup of full shift space does not form a sensitive system (under the action of shift operator).
\end{cor}

\begin{proof}
As the full shift is not weakly rigid (contains non-recurrent points), the enveloping semigroup fails to form a sensitive flow (by corollary \ref{wr}).
\end{proof}

\begin{remark}
The above results relate the sensitivity of the flow $(E(X), T)$ with weak rigidity of the underlying flow $(X, T)$. The proof utilizes the fact that a sensitive flow cannot contain any isolated element in it. In particular, we can conclude that the enveloping semigroup of a sensitive flow may not form a sensitive flow. It is known that any homomorphism between minimal flows is feeble open~\cite{a}. It is known that for compact metrizable flows, if any factor of $(X,T)$ under a feeble open factor map is sensitive then the flow $(X, T)$ is also sensitive \cite[Proposition 7.2.4.]{Vries}. We establish the result under a more general setting of compact Hausdorff spaces (Theorem \ref{sensitive}). Consequently, if $(X, T)$ is distal and has a minimal sensitive subsystem then the induced flow $(E(X),T)$ is sensitive (Theorem \ref{sub}). However, it can be seen that sensitivity of the induced flow $(E(X),T)$ in general does not guarantee sensitivity for $(X,T)$. 
\end{remark}

\begin{theorem}\label{sensitive}
Let $\phi:(X, T)\rightarrow (Y, T)$ be a factor map of flows with $(X, T)$ and $(Y, T)$ compact Hausdorff spaces (not necessarily metrizable) and assume that $\phi$ is feeble-open. If $(Y, T)$ is sensitive, then $(X, T)$ is also sensitive.
\end{theorem}

\begin{proof}
Let $\mathcal{U}_X$ and $\mathcal{U}_Y$ be uniform structure on $X$ and $Y$ respectively.
Since $(Y,T)$ is sensitive, there exists $V\in \mathcal{U}_Y$ such that for any $y\in Y$ and any neighbourhood $N_y$ of $y$, there exists $y^{'}\in N_y$ and $g\in T$ such that $(gy,gy^{'}) \notin V$. Further, as $X$ is a compact space, $\phi$ is uniformly continuous. Consequently there exists $U\in \mathcal{U}_X$ such that for all $(x,x^{' })\in U$, $(\phi(x),\phi(x^{'}))\in V$. Without loss of generality, let $U=U_0\circ U_0$. Now, we will show that $(X, T)$ is $U_0-sensitive$. If not, then there exists $x_0\in X$ and a neighborhood $N_{x_0}$ of $x_0$ such that $(gx,gx_{0})\in U_0$ for all $g\in T$ and $x\in N_{x_0}$ which further implies $(gx,gx^{'})\in U$, for all $g\in T$ and $x,x^{'}\in N_{x_0}$. Since $\phi$ is semi-open, $\phi(N_{x_0})$ includes a non-empty open set, say $W$. Then for $y,y^{'}\in W$, choose $x,x^{'}\in N_{x_0}$ such that $\phi(x)=y$ and $\phi(x^{'})=y^{'}$. Then for $g\in T$, we have $(gy,gy^{'})=(g\phi(x),g\phi(x^{'}))=(\phi(gx),\phi(gx^{'}))\in V$, which is a contradiction to the fact that $(Y,T)$ is sensitive and hence the proof is complete.
\end{proof}

\begin{theorem}\label{sub}
If $(X, T)$ is distal and has a minimal sensitive subsystem, then $(E(X), T)$ is sensitive.
\end{theorem}
\begin{proof}
Let $Y=\overline{\mathcal{O}(y)}$ be minimal sensitive subsystem of $(X,T)$. Define $\phi:(E(X),T)\rightarrow \overline{\mathcal{O}(y)}$ as $\phi(p)=py$. It may be noted that $\phi$ is a surjective homomorphism. Also, as the flow $(X, T)$ is distal, $(E(X), T)$ is a minimal flow and hence $\phi$ is a homomorphism of minimal flows, which further implies $\phi$ is a feeble-open map. Since $Y$ is a sensitive flow, we can conclude that $(E(X), T)$ forms a sensitive flow (Theorem~\ref{sensitive}).
\end{proof}

\begin{example}
Let $X=\mathbb{T}\cup (S^1\times\{0\})\subset \R^3$, where $\mathbb{T}$ is a 2-dimensional torus (or $S^1\times S^1$) disjoint from $S^1\times\{0\}$ (represented as in Figure ~\ref{fig:4}) . Now define $f:X\rightarrow X$ as $f|_{\mathbb{T}}(\theta_1,\theta_2)=(\theta_1+\mu,\theta_1+\theta_2)$ and $f|_{S^1\times\{0\}}(\theta,0)=(\theta+\alpha,0)$, where $\alpha$ and $\mu$ are not roots of unity.
\begin{figure}[hbt!]
    \centering
    \includegraphics[width=0.9\linewidth]{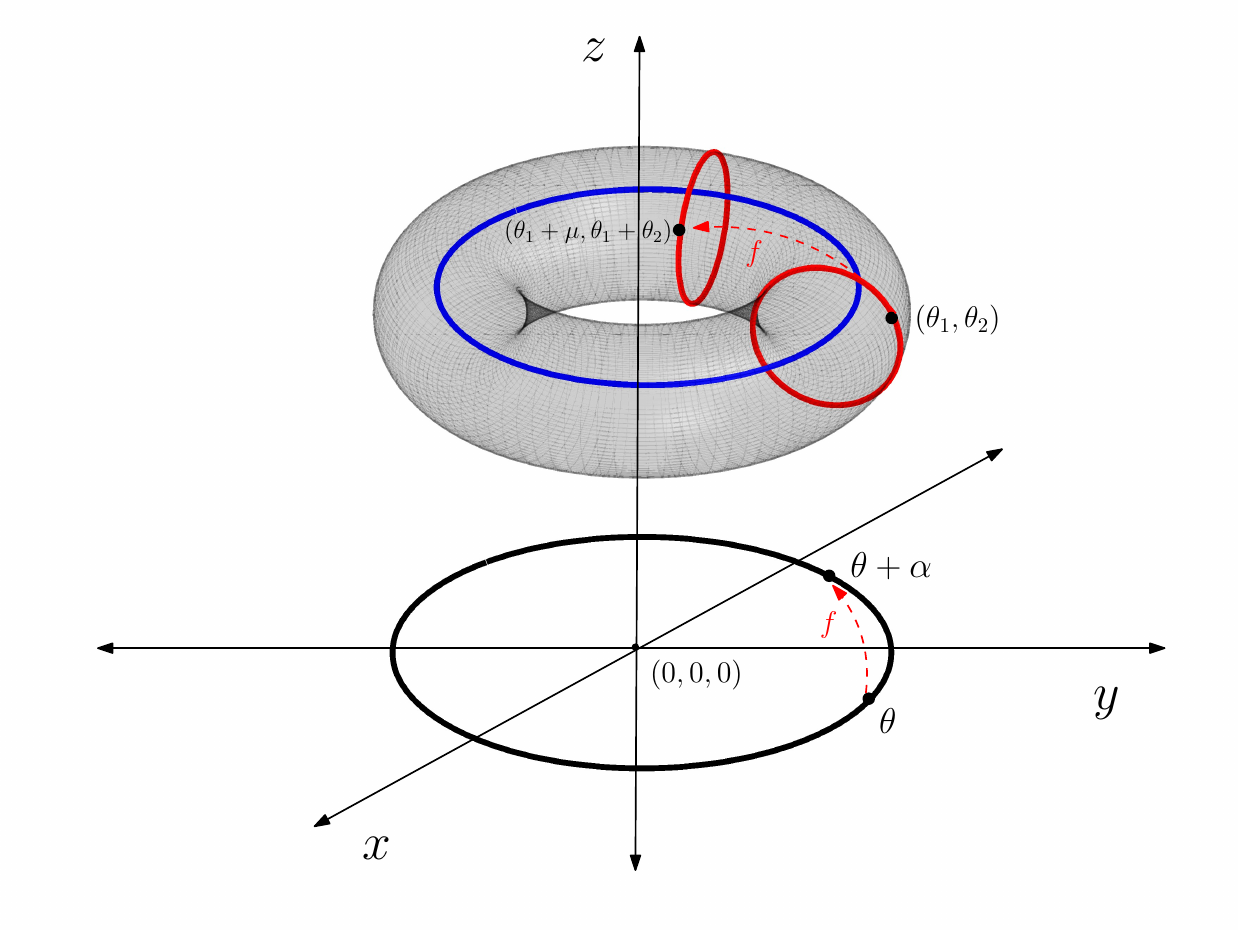}
    \caption{}
    \label{fig:4}
\end{figure}
Note that the system $(X,f)$ is distal, and it has a minimal sensitive subsystem $(\mathbb{T},f)$. Let $\delta>0$ be sensitivity constant for $(\mathbb{T},f)$. Now we claim that for any fixed $y\in \mathbb{T}$, $V=\{(f,g)\in E(X)\times E(X): d(f(y),g(y))<\frac{\delta}{2}\}$ is sensitive entourage for $(E(X),\hat{f})$, where $\hat{f}(g)=f\circ g$, for $g\in E(X)$. Let $g\in E(X)$ be any element and $S_g(x,\epsilon)=\{h\in E(X): d(g(x),h(x))<\epsilon\}$ be any sub-basic neighborhood of $g$. For any fixed $y\in \mathbb{T}$, as $\phi: E(X)\rightarrow \mathbb{T}$ given by $\phi(g)=gy$ is a surjective homomorphism between minimal systems, $\phi$ is feeble open~\cite{a} and thus $S_g(x,\epsilon)(y)=\{h(y): h\in S_g(x,\epsilon)\}$ has non empty interior (say $W$) in $\mathbb{T}$. Since $(\mathbb{T},f|_{\mathbb{T}})$ is sensitive, there exists $y_1,y_2\in W$ and $k\in \N$ such that $d(f^k(y_1),f^k(y_2))>\delta$. Consequently there exists $h_1,h_2\in S_g(x,\epsilon)$ such that $h_1(y)=y_1$, $h_2(y)=y_2$ and $d(f^k(h_1(y)),f^k(h_2(y)))>\delta$. Also note that either $(\hat{f^k}(g),\hat{f^k}(h_1))\notin V$ or $(\hat{f^k}(g),\hat{f^k}(h_2))\notin V$ and hence $(E(X),\hat{f})$ is sensitive with sensitive entourage $V$.

\end{example}

\begin{remark}
It is known that for compact metrizable flows, any point-transitive system is either almost equicontinuous or sensitive. Also, it is known that any point transitive almost equicontinuous is uniformly rigid~\cite{AB}. As compact Hausdorff spaces have a unique uniformity generating the undelying topology, the results hold under a more general setting of compact Hausdorff spaces (Theorem \ref{ad2} and Theorem \ref{ur}). We now establish our claims below.
\end{remark}

\begin{theorem}\label{ad2}
If $(X, T)$ is a point-transitive flow with no isolated points, then $(X,T)$ is either almost equicontinuous or sensitive. Further, if $(X,T)$ is almost equicontinuous, then $Eq(X)=Trans(X)$.
\end{theorem}

\begin{proof}
Let $(X, T)$ be a point-transitive system with no isolated points and $x_0\in X$ be a transitive point. Firstly, note that if $x_0$ is a point of equicontinuity, then $gx_0$ is also a point of equicontinuity (for any $g\in T$), and hence $(X,T)$ is almost equicontinuous. Further, if $x_0$ is a point of sensitivity with sensitive entourage $V\in \mathcal{U}_X$, then $gx_0$ is also a point of sensitivity with sensitive entourage $V\in \mathcal{U}_X$. The proof follows from the fact that for any entourage $U\in \mathcal{U}_X$ there exists $y\in g^{-1}U[x_0]$ and $t\in T$ such that $(tx_0,ty)\notin V$. Consequently, $(tg^{-1}gx_0,tg^{-1}gy)\notin V$ and thus $gx_0$ is also a point of sensitivity. Let $z\in X$ be fixed and $U\in \mathcal{U}_X$ be any entourage. Since $x_0$ is a transitive point, there exists $g\in T$ such that $gx_0\in U[z]$. Also, there exists $y\in U[g(x_0)]$ and $t\in T$ such that $(tgx_0,ty)\notin V$. Further, if $V_0\in \mathcal{U}_X$ is such that $V_0\circ V_0\subseteq V$, then either $(tgx, tz)\notin V_0$ or $(tz, ty)\notin V_0$ holds and thus $z$ is a point of sensitivity with sensitive entourage $V_0$. As the proof holds for any $z\in X$, $(X,T)$ is sensitive. \\ 

Finally, if $(X, T)$ is almost equicontinuous, then from the above discussion, we can conclude that $Trans(X)\subseteq Eq(X)$. Let $x_0\in X$ be a point of equicontinuity, $z\in X$ be any element and $V\in \mathcal{U}_X$ be any entourage. For $V\in \mathcal{U}_X$, choose a symmetric entourage $V_0\in \mathcal{U}_X$ such that $V_0\circ V_0\subseteq V$. Since $x_0\in X$ is a point of equicontinuity, there exists $U\in \mathcal{U}_X$ such that $(x_0,y)\in U$ ensures $(tx_0,ty)\in V_0$ for all $t\in T$ and $y\in X$. Also, if  $y_0\in X$ is a transitive point (as $(X,T)$ is point-transitive), then there exists $t_0\in T$ such that $t_0y_0\in U[x_0]$ and consequently $(tt_0y_0,tx_0)\in V_0$ for all $t\in T$. Consequently, if $t_1\in T$ is such that $t_1y_0\in V_0[z]$ then $(t_1t_0^{-1}x_0,z)\in V$. Hence, $x_0$ is a transitive point, and we have $Eq(X)=Trans(X)$. 
\end{proof}

\begin{theorem}\label{ur}
If $(X, T)$ is point-transitive, almost equicontinuous flow such that $T$ is abelian, then $(X, T)$ is uniformly rigid. 
\end{theorem}
\begin{proof}
Let $(X, T)$ be a point-transitive almost equicontinuous flow, $x_0\in X$ be a point of equicontinuity and $V\in \mathcal{U}_X$ be any entourage. Choose $V_0\in \mathcal{U}_X$ such that $\overline{V_0}\subset V$. Since $x_0$ is a point of equicontinuity, for $V_0\in \mathcal{U}_X$, there exists $U\in \mathcal{U}_X$ such that $(tx_0,ty)\in V_0$ for all $y\in U[x_0]$ and $t\in T$. Also as $x_0$ is transitive (Theorem \ref{ad2}), there exists $t_0\in T$ such that $t_0x_0\in U[x_0]$ and thus $(tx_0,t_0tx_0)\in V_0$ for all $t\in T$ (as $T$ is abelian). Further for any $y\in X$, there exists a net $t_{\alpha}$ in $T$ such that $t_\alpha x_0\rightarrow y$. As $(t_\alpha x_0,t_0t_\alpha x_0)\in V_0$ (for all $\alpha$), we have $(y,t_0 y)\in \overline{V_0}\subset V$. As the proof holds for any $y\in X$, $(X,T)$ is uniformly rigid and the proof is complete.
\end{proof}

\begin{cor}
 For any point-transitive, weakly rigid abelian flow $(X,T)$ which is not uniformly rigid, $(E(X), T)$ forms a sensitive flow.
\end{cor}

\begin{proof}
Firstly note that as $(X, T)$ is weakly rigid, $E(X)$ has no isolated point (Theorem \ref{niso}) and thus is either sensitive or almost equicontinuous (Theorem~\ref{ad2}). Also note that in case $(E(X), T)$ is almost equicontinuous, then it is uniformly rigid (Theorem~\ref{ur}), which forces $(X, T)$ is uniformly rigid (as $(X, T)$ is a factor of $(E(X), T))$ which is a contradiction. Hence, $(E(X), T)$ forms a sensitive flow and the proof is complete.
\end{proof}

\begin{remark}
Note that the flow $(X, T)$, which is point-transitive, weakly rigid but not uniformly rigid is a sensitive flow. Thus the above result provides a sufficient condition for the enveloping semigroup of a sensitive flow to be sensitive. Further, it can be seen that the above result fails to hold in the absence of point transitivity (Example \ref{circ}).
\end{remark}

\begin{theorem}\label{weakly}
Any flow $(X, T)$ is weakly rigid if and only if $(E(X), T)$ is weakly rigid.
\end{theorem}

\begin{proof}
Let $(X, T)$ be a weakly rigid flow and $(p_1,p_2,\ldots,p_n)\in E(X)^n$. Also, consider $\Pi_{i=1}^nU(x_i,V_i)$ be sub-basic open set containing $(p_1,p_2,\ldots,p_n)$ (where $x_i\in X$ and $V_i$ be open set in $X$ for all $i\in\{1,2,\ldots,n\}$). Since $(X, T)$ is a weakly rigid flow, there exists $t\in T$ such that $tp_i(x_i)\in V_i$ and hence $(tp_1,tp_2,\ldots,tp_n)\in \Pi_{i=1}^n U(x_i,V_i)$. Consequently, $(E(X), T)$ is weakly rigid. Conversely, consider $x_1,x_2,\ldots,x_n\in X$ and $\Pi_{i=1}^nU_i$ be open set in $X^n$ containing $(x_1,x_2,\ldots,x_n)$. As $\Pi_{i=1}^nU(x_i,U_i)$ is a neighborhood of $(e,e,\ldots,e)\in E(X)^n$, there exists $t\in T$ such that $te(x_i)\in U_i$ (as $(E(X),T)$ is weakly rigid) and hence $(X,T)$ is weakly rigid.
\end{proof}
\begin{theorem}\label{uni}
If a flow $(X, T)$ is uniformly rigid, then $(E(X), T)$ also forms a uniformly rigid flow. Further, if $(X, T)$ is point-transitive, the converse also holds.
\end{theorem}

 \begin{proof}
Let $S(x_0,U)=\{(f,g):(f(x_0),g(x_0))\in U)\}$ (where $x_0\in X$ and $U\in \mathcal{U}_X$ are fixed) be sub-basic element of uniformity on $E(X)$. Since $(X, T)$ is uniformly rigid, there exists a net $(\pi^{t_\alpha})$ converging uniformly to $e$. Consequently, there exists $\beta\in J$(directed set) such that $(\pi^{t_\alpha}(x),x)\in U$ for all $\alpha\geq \beta$, $x\in X$ which implies $(\pi^{t_\alpha}(g(x_0)),g(x_0))\in U$ for all $g\in E(X)$ and $\alpha\geq \beta$. Thus, $(\hat{\pi}^{t_\alpha}(g),g)\in S(x_0,U)$ for all $g\in E(X)$ and hence $(E(X), T)$ forms a uniformly rigid flow. Conversely, if $(X, T)$ is point-transitive, then $(X, T)$ is a factor of $(E(X), T)$ and hence is uniformly rigid (as a factor of uniformly rigid flow is uniformly rigid).
\end{proof}

\begin{remark}
The above results establish that if the flow $(X, T)$ exhibits weak (uniform) rigidity, then $(E(X),T)$ also exhibits the same. However, if  $(E(X),T)$ is uniformly rigid, $(X,T)$ need not exhibit the same. For instance, in Example~\ref{circ}, $(X,f)$ is distal and hence weakly rigid. Also, as $(X,f)$ is hereditarily almost equicontinuous, $(E(X),f)$ is equicontinuous and hence uniformly rigid (Theorem \ref{t3}).  However as $f^{2^n}$ converges pointwise to $e$ (not uniformly) on $X$,  $(X,f)$ is not uniformly rigid. Further, it is known that any compact metrizable system $(X,f)$ is proximal if and only if it is syndetically proximal \cite[Theorem 1]{SM}. A similar result can be established for a general system $(X, T)$ (where $X$ need not be metrizable and $T$ is abelian). %We include the proof for the sake of completion.
\end{remark}

\begin{theorem}\label{synd}
For any flow $(X,T)$ where $T$ is abelian, following are equivalent:
\begin{enumerate}
\item $(X,T)$ is proximal.
\item $(X, T)$ has a fixed point as a unique minimal subset.
\item $(X,T)$ is syndetically proximal.
\end{enumerate}
\end{theorem}

\begin{proof}
Firstly note that if $(X,T)$ is proximal, then for any $x_0\in X$ and $t_0\in T$, the pair $(x_0,t_0x_0)$ is proximal and hence there exists a net $(t_\alpha)$ such that $\lim t_\alpha x_0=\lim t_\alpha t_0x_0 = z_0$(say). As $T$ is abelian, $z_0$ is a fixed point under the action of $t_0$. Thus for any $t^*\in T$, the set $A_{t^*}=\{x\in X: t^*x=x\}$ is non-empty, closed invariant subset of $X$ (as $T$ is abelian). As $(X, T)$ is proximal, any two closed invariant sets (and hence finitely many closed invariant sets) must intersect. Consequently, $\{A_t:t\in T\}$  has finite intersection property and thus $\bigcap \limits_{t\in T} A_t \neq\emptyset$ (as $X$ is compact). Thus there exists a point $p\in X$ such that $tp=p$ for all $t\in T$. As the existence of two fixed points (minimal subsystems) violates the proximality of $(X, T)$, $(X, T)$ has a fixed point as a unique minimal subset and $(1)\implies (2)$ holds.\\

If the flow  $(X, T)$ has a fixed point, say $z\in X$ as a unique minimal subset, then $(z,z)$ forms a unique minimal set in $X\times X$. As $(z,z)\in \overline{\mathcal{O}(x,y)}$ (for $(x,y)\in X\times X$), for any $U\in\mathcal{U}_X$ there exists $t\in T$ such that $(tx,ty)\in U$ and thus we have $\cup_{t\in T}tU=X\times X$. As $X$ is compact, there exists $K=\{t_1,t_2,\ldots,t_n\}$ such that $\cup_{t\in K} tU=X\times X$. Hence for any $(x,y)\in X\times X$ and $U\in \mathcal{U}_X$, $A=\{t:(tx,ty)\in U\}$ is syndetic (as $KA=T$)  and $(X,T)$ is syndetically proximal. Hence $(2)\implies (3)$ holds.\\

As $(3)\implies (1)$ holds trivially, the proof is complete.
\end{proof}

\begin{theorem}\label{ts}
For a flow $(X, T)$, any pair $(x,y)\in X\times X$ is a syndetically proximal pair if and only if for any $V\in \mathcal{U}_X$, $\{t\in T:(tx, ty)\in V\}$ is a thickly syndetic set.   
\end{theorem}
\begin{proof}
Let $(x,y)$ be a syndetically proximal pair, $V\in \mathcal{U}_X$ be an entourage and $K\subset T$ be any compact set. %Since $(x,y)$ is syndetically proximal, $A=\{t\in T:(tx, ty)\in V\}$ is a syndetic set. 
As $K$ is compact, there exists $U\in\mathcal{U}_X$ such that $(x,y)\in U$, $k\in K$ implies $(kx,ky)\in V$. Also as $(x,y)$ is syndetically proximal, $S_U=\{t\in T:(tx, ty)\in U\}$ is a syndetic set. Consequently, if $S_V=\{t\in T:(tx, ty)\in V\}$ then $KS_U\subseteq S_V$ and hence $S_V$ is thickly syndetic. \\

Conversely, for any pair $(x,y)\in X\times X$ and $V\in \mathcal{U}_X$, if $\{t\in T:(tx, ty)\in V\}$ is thickly syndetic then  $(x,y)\in X\times X$ is a syndetically proximal and the proof is complete.
\end{proof}

\begin{theorem}
For abelian $T$, if $(X, T)$ is a proximal flow, then $(E(X), T)$ is also a proximal flow. Further, if $(X, T)$ is point-transitive, then the converse is also true.
\end{theorem}

\begin{proof}
Let $(f,g)$ be any pair in $E(X)\times E(X)$ and $V=\cap_{i=1}^{n}S(x_i,U_i)\in \mathcal{U}_{E(X)}$ be any basis element in uniformity of $E(X)$ (where $x_i\in X$ and $U_i\in \mathcal{U}_X$ for all $i=1,2,\ldots, n$). As $(X,T)$ is proximal (hence syndetically proximal), for each pair $(f(x_i),g(x_i))$ the set $ A_i=\{t\in T:(tf(x_i),tg(x_i))\in U_i\}$ is thickly syndetic (Theorem~\ref{synd} and Theorem~\ref{ts}). Further, we know that the intersection of thickly syndetic sets is again a thickly syndetic set and thus we have the set $\{t\in T: (tf,tg)\in V\}$ is a thickly syndetic set and hence $(E(X), T)$ is a proximal (syndetically proximal) flow. Further, if $(X, T)$ is point-transitive, then $(X,T)$ is a factor of $(E(X),T)$, and hence $(X, T)$ is proximal (as a factor of any proximal flow is proximal).
\end{proof}

\begin{remark}
The above result relates proximality of flow $(X, T)$ and $E(X, T)$ when the acting group $T$ is abelian. The result establishes that if a flow $(X, T)$ is point-transitive, then proximality of $(E(X),T)$ implies proximality of $(X,T)$.  However, the conclusion does not hold in the absence of point transitivity of $(X,T)$. We now give an example in support of our claim.
\end{remark}

\begin{example}

Let $X=\{0,1\}^\Z$ be the collection of all bi-infinite sequences over two symbols and $\sigma$ be the left shift operator on $X$. Consider $b=x_1x_2\ldots x_n$ be any fixed block, $x=\{0^\infty .b 0^\infty\}$ and $y=\{1^\infty .b 1^\infty\}$. Then if $Y=\overline{\mathcal{O}(x)}\cup \overline{\mathcal{O}(y)}$ then $\Z$-action defined on $Y$ as $1.x=\sigma(x)$ forms a flow under action of left shift operator. Also $(0^{\infty},1^{\infty})\in Y\times 
Y$ is not proximal for $(Y,T)$ and hence $(Y,T)$ is not a proximal flow. However, $E(Y)=\{\sigma^n:n\in \Z\}\cup\{g\}$ where $g:Y\rightarrow Y$ is defined as $$g(s) = \left\{
     \begin{array}{lr}
       {0^{\infty}} &  : s\in \overline{\mathcal{O}(x)}  \\
       {1^{\infty}} & : s\in \overline{\mathcal{O}(y)} \\
       \end{array}\right\}$$

As $\sigma^n\rightarrow g$ as $n\rightarrow \pm\infty$, $(E(Y),T)$ forms a proximal flow.
\end{example}

\begin{theorem}
If $(X,T)$ is syndetically distal, $(E(X),T)$ is syndetically distal.
\end{theorem}

\begin{proof}
Let $f,g\in E(X)$  be distinct in $E(X)$ such that $(f,g)$ forms a syndetically proximal pair. As $f$ and $g$ are distinct, there exists $x\in X$ such that $f(x)\neq g(x)$. Also, for any $S(x,U)\in \mathcal{U}_{E(X)}$ (where $U\in \mathcal{U}_X$), there exists a syndetic set $A\subseteq T$ such that $(tf,tg)\in S(x,U)$ for all $t\in A$ which implies $(tf(x),tg(x))\in U$ for all $t\in A$. As the proof holds for any $U\in \mathcal{U}_X$, $(f(x),g(x))$ is a syndetically proximal pair in $(X, T)$, which contradicts the syndetic distality of $(X,T)$.
\end{proof}

\begin{remark}\label{last}
The above proof establishes that if $(X,T)$ is syndetically distal, then $(E(X),T)$ is syndetically distal. Conversely, if $(E(X),T)$ is syndetically distal, then $(E(X),T)$ contains a distal pair, and hence $(X,T)$ contains a distal pair. However, as distality of $(X,T)$ implies distality of $(E(X),T)$ if $(E(X),T)$ is syndetically distal but not distal, then $(X,T)$ contains both non-trivial proximal and distal pairs. Thus we get the following corollary.
\end{remark}

\begin{cor}
If $(E(X),T)$ is syndetically distal, then $(X,T)$ is not proximal (but contains non-trivial proximal pairs).
\end{cor}

\begin{proof}
The proof follows from discussions in Remark \ref{last}.
\end{proof}

\section*{Acknowledgement}
The first author thanks MoE (India) for financial support. The second author thanks the National Board for Higher Mathematics (NBHM) Grant No. 02011/27/2023NBHM(R.P)/R \& D II/6282 for financial support.

\end{document}